\documentclass[a4paper,11pt,twoside,english]{article}
\usepackage{babel}
\usepackage{latexsym,amsfonts}
\usepackage{amsthm}
\usepackage{amsmath}
\usepackage{mathtools}
\usepackage{paralist}

\usepackage{lastpage}
\usepackage{fancyhdr}
\fancypagestyle{plain}{

\fancyhead[L,C]{}
\fancyhead[R]{\today}
\fancyfoot[R,L]{}
\fancyfoot[C]{\thepage \ of \pageref{LastPage}}
}

\title{Some Qi-type integral inequalities involving several weight functions}
\author{Jan-David Hardtke}
\date{}

\providecommand{\sm}{\setminus}

\providecommand{\N}{\ensuremath{\mathbb{N}}}

\providecommand{\R}{\ensuremath{\mathbb{R}}}

\providecommand{\keywords}[1]{
{\let\thefootnote=\relax
\footnote{{\em Keywords}: #1}}
\addtocounter{footnote}{-1}
}

\providecommand{\AMS}[1]{
{\let\thefootnote=\relax
\footnote{{\em AMS Subject Classification} (2010): #1}}
\addtocounter{footnote}{-1}
}

\providecommand{\address}{
{\sc \noindent Department of Mathematics \\
Freie Universit\"at Berlin \\
Arnimallee 6, 14195 Berlin \\
Germany \\}
}

\DeclarePairedDelimiter{\set}{\lbrace}{\rbrace}
\DeclarePairedDelimiter{\paren}{\lparen}{\rparen}

\DeclarePairedDelimiter{\norm}{\lVert}{\rVert}

\theoremstyle{definition}
\newtheorem{definition}{Definition}
\newtheorem*{definition*}{Definition}
\theoremstyle{plain}
\newtheorem{lemma}[definition]{Lemma}
\newtheorem*{lemma*}{Lemma}
\newtheorem{proposition}[definition]{Proposition}
\newtheorem*{proposition*}{Proposition}

\newtheorem*{theorem*}{Theorem}
\newtheorem{corollary}[definition]{Corollary}
\newtheorem*{corolary*}{Corollary}

\newenvironment{Proof}[1][\proofname]{\begin{proof}[#1] \setlength{\parindent}{0pt}}{\end{proof}}
\newenvironment{Abstract}{\centering\begin{minipage}{0.8\textwidth} \noindent \small {\sc Abstract.}}{\end{minipage}\par}

\usepackage{color}
\definecolor{darkgreen}{rgb}{0,0.5,0}

\newtagform{colored}[\color{blue}]{\color{blue}(}{\color{blue})}
\usetagform{colored}

\usepackage[colorlinks,linkcolor=blue,citecolor=red,urlcolor=darkgreen]{hyperref}
\providecommand{\email}{{\it E-mail address:} \href{mailto:hardtke@math.fu-berlin.de}{\tt hardtke@math.fu-berlin.de}}

\usepackage{amsrefs}

\begin{document}

\maketitle

\begin{Abstract}
We prove some integral inequalities related to Feng Qi's inequality from \cite{qi} and obtain a few corollaries.
\end{Abstract}
\AMS{26D15}

In the paper \cite{qi}, Feng Qi proved the following integral inequality: if $f:[a,b] \rightarrow \R$ is $n$-times continuously
differentiable such that $f^{(i)}(a)\geq 0$ for $i=0,\dots,n-1$ and $f^{(n)}(x)\geq n!$ for all $x\in [a,b]$, then
\begin{equation}\label{eq:qi}
\paren*{\int_a^b f}^{n+1}\leq \int_a^b f^{n+2}.
\end{equation}
In the following years, many variants and generalisations of this inequality have been investigated (in particular, versions for
real exponents). For a detailed account of the results that have already been established, the reader is referred to the list of 
references at the end of this paper and further references therein.\par
Here we will prove some further inequalities related to Feng Qi's inequality and obtain a few corollaries which might be interesting.
The basic methods of proof are the same as in \cite{qi} (exploiting the connection between a function's monotonicity behavior and the 
sign of its derivative, finite induction for a suitably defined set of auxiliary functions, etc.).\par
Since we will also be dealing with one-sided derivatives, let us introduce some notations and recall some facts, which will 
be used later without further mention.\par 
For a given interval $[a,b]$ we denote by $D_{+}[a,b]$ (resp. $D_{-}[a,b]$) the set of all continuous functions 
$h:[a,b] \rightarrow \R$ such that the right-derivative $h_{+}^{\prime}(x)$ (resp. the left-derivative $h_{-}^{\prime}(x)$) 
exists for all $x\in (a,b)$.\par 
The usual sum and product rules also apply for one-sided derivatives. Furthermore, 
a function $h\in D_{+}[a,b]$ is increasing (decreasing) if and only if $h_{+}^{\prime}(x)\geq 0$ ($h_{+}^{\prime}(x)\leq 0$) 
for all $x\in (a,b)$ and an analogous statement holds for functions in $D_{-}[a,b]$ (see \cite{walter}*{p.358}).\par
Finally, if $h\in D_{+}[a,b]$, $I$ is an interval which contains the range of $h$ and
$f:I \rightarrow \R$ is differentiable, then $f\circ h\in D_{+}[a,b]$ with
$(f\circ h)_{+}^{\prime}(x)=f^{\prime}(h(x))h_{+}^{\prime}(x)$ for all $x\in (a,b)$
(and analogously for functions in $D_{-}[a,b]$). This is proved in the same manner
as the usual chain-rule.\par 

\section{A generalisation of Qi's inequality}\label{sec:general-qi}
	
The first result is a generalisation of \eqref{eq:qi} to a setting where several weight functions and an additional exponent $\alpha$ are involved.

\begin{proposition}\label{prop:general-qi-inequal}
Let $n\in \N$ with $n\geq 2$ and let $g:[a,b] \rightarrow \R$ be a continuous function with $g\geq 0$. Further, let $h_1,\dots,h_n\in D_{+}[a,b]$ or 
$h_1,\dots,h_n\in D_{-}[a,b]$ with $h_i\geq 0$ for each $i$. Let also 
$\alpha\leq \frac{n}{n-1}$ and $f:[a,b] \rightarrow \R$ be a strictly positive, 
$(n-1)$-times differentiable function such that $f^{(i)}\geq 0$ for all $i=1,\dots,n-1$ and $(n+1-\alpha)f^{(n-1)}f^{n(1-\alpha)}\geq n!h_1h_2^2\dots h_n^n$.\par
Suppose that there exists a partition of $\set*{1,\dots,n}$ into two disjoint subsets $I$ and $J$ such that the following conditions hold:
\begin{enumerate}[\upshape(i)]
\item $g\prod_{i\in I}h_i$ is increasing,
\item $\paren*{\prod_{i\in I}h_i^i}\paren*{\prod_{j\in \underline{J}_k}h_j^{j-k}}$ is decreasing for all $k=1,\dots,n-1$, where $\underline{J}_k$ denotes the set 
$\set*{j\in J:j\geq k+1}$.
\item $\prod_{j\in \overline{J}_k}h_j^{k+1-j}$ is increasing for $k=0,\dots,n-2$, where $\overline{J}_k$ denotes the set $\set*{j\in J:j\leq k+1}$.
\end{enumerate}
Let $M:=\inf\set*{g(x):x\in [a,b]}$ and $m_i:=\inf\set*{h_i(x):x\in [a,b]}$ for $i=1,\dots,n$.
Then we have
\begin{equation*}
(b-a)KM^n\prod_{i=1}^nm_i^{n-i}+\paren*{\int_a^bf^{\alpha}g\prod_{i=1}^nh_i}^n\leq\int_a^bf^{n+1}g^n\prod_{i=1}^nh_i^{n-i},
\end{equation*}
where $K=f^{n+1}(a)$ if $\alpha\geq 0$ and $K=f^{n+1-\alpha}(a)f^{\alpha}(b)$ if $\alpha<0$.
\end{proposition}

Note that the conditions (ii) and (iii) are satisfied in particular if $h_1h_2^2\dots h_n^n$ is decreasing and $h_j$ is increasing and strictly 
positive for each $j\in J$.\par
Note further that the above inequality not only generalises Qi's inequality, but it also sharpens the trivial estimate
\begin{equation*}
(b-a)KM^n\prod_{i=1}^nm_i^{n-i}\leq\int_a^bf^{n+1}g^n\prod_{i=1}^nh_i^{n-i}.
\end{equation*}

\begin{Proof}
Let $h_1,\dots,h_n\in D_{+}[a,b]$ (the case of left-derivatives is treated analogously).\par 
We first define the function $F$ by setting 
\begin{equation*}
F(x):=\paren*{\int_a^xf^{\alpha}g\prod_{i=1}^nh_i}^n-\int_a^xf^{n+1}g^n\prod_{i=1}^nh_i^{n-i} \ \ \ \forall x\in [a,b].
\end{equation*}
Then $F$ is differentiable with
\begin{equation}\label{Eq:diffF}
F^{\prime}(x)=n\paren*{\int_a^xf^{\alpha}g\prod_{i=1}^nh_i}^{n-1}f^{\alpha}(x)g(x)\prod_{i=1}^nh_i(x)-f^{n+1}(x)g^n(x)\prod_{i=1}^nh_i^{n-i}(x).
\end{equation}
Next we define
\begin{equation*}
G(x):=n\paren*{\int_a^xf^{\alpha}\prod_{j\in J}h_j}^{n-1}\paren*{\prod_{i\in I}h_i^i(x)}\paren*{\prod_{j\in J}h_j(x)}-f^{n+1-\alpha}(x)\prod_{j\in J}h_j^{n-j}(x)
\end{equation*}
and claim that 
\begin{equation}\label{Eq:FG}
F^{\prime}(x)\leq f^{\alpha}(x)g^n(x)G(x)\prod_{i\in I}h_i^{n-i}(x) \ \ \forall x\in [a,b].
\end{equation}
To see this, note that assumption (i) implies
\begin{equation*}
\int_a^xf^{\alpha}g\prod_{i=1}^nh_i\leq g(x)\paren*{\prod_{i\in I}h_i(x)}\int_a^xf^{\alpha}\prod_{j\in J}h_j.
\end{equation*}
Combining this with \eqref{Eq:diffF} we obtain
\begin{align*}
&F^{\prime}(x)\leq ng^n(x)f^{\alpha}(x)\paren*{\int_a^xf^{\alpha}\prod_{j\in J}h_j}^{n-1}\paren*{\prod_{i\in I}h_i^{n-1}(x)}\paren*{\prod_{i=1}^nh_i(x)} \\
&-f^{n+1}(x)g^n(x)\prod_{i=1}^nh_i^{n-i}(x)
\end{align*}
which can be easily simplified to \eqref{Eq:FG}.\par
We denote by $\varphi$ the characteristic function of $J$ in $\set*{1,\dots,n}$, i.\,e. $\varphi(i)=1$ for $i\in J$ and $\varphi(i)=0$ for $i\in I$,
and for each $k\in \set*{1,\dots,n-1}$ we define
\begin{align*}
&H_k(x):=\paren*{\prod_{i=k}^ni}\paren*{\int_a^xf^{\alpha}\prod_{j\in J}h_j}^{k-1}\paren*{\prod_{i\in I}h_i^i(x)}\paren*{\prod_{i=k}^nh_i^{\varphi(i)(i-k+1)}(x)} \\
&-(n+1-\alpha)f^{n(1-\alpha)+(k-1)\alpha}(x)f^{(n-k)}(x)\prod_{i=1}^{k-1}h_i^{(k-i-1)\varphi(i)}(x) \ \ \ \text{for\ all} \ x\in [a,b].
\end{align*}
We will show inductively that $H_k\leq 0$ for all $k\in \set*{1,\dots,n-1}$. For $k=1$ we have $H_1=n!\prod_{i=1}^nh_i^i-(n+1-\alpha)f^{(n-1)}f^{n(1-\alpha)}$,
which is negative by assumption.\par
Now suppose that $k\in\set*{1,\dots,n-2}$ and $H_k\leq 0$. The function $H_{k+1}$
belongs to $D_{+}[a,b]$ and satisfies
\begin{align*}
&(H_{k+1})_{+}^{\prime}(x)= \\
&f^{\alpha}(x)\paren*{\prod_{i=k}^ni}\paren*{\int_a^xf^{\alpha}\prod_{j\in J}h_j}^{k-1}\paren*{\prod_{j\in J}h_j(x)}\paren*{\prod_{i\in I}h_i^i(x)}
\paren*{\prod_{i=k+1}^nh_i^{\varphi(i)(i-k)}(x)} \\
&+\paren*{\prod_{i=k+1}^ni}\paren*{\int_a^xf^{\alpha}\prod_{j\in J}h_j}^k\paren*{\paren*{\prod_{i\in I}h_i^i}\paren*{\prod_{i=k+1}^nh_i^{\varphi(i)(i-k)}}}_{+}^{\prime}(x) \\
&-(n+1-\alpha)\bigg((n(1-\alpha)+k\alpha)f^{n(1-\alpha)+k\alpha-1}(x)f^{\prime}(x)f^{(n-k-1)}(x)\prod_{i=1}^{k}h_i^{\varphi(i)(k-i)}(x) \\
&+f^{n(1-\alpha)+k\alpha}(x)f^{(n-k)}(x)\prod_{i=1}^{k}h_i^{\varphi(i)(k-i)}(x) \\
&+f^{n(1-\alpha)+k\alpha}(x)f^{(n-k-1)}(x)\paren*{\prod_{i=1}^{k}h_i^{\varphi(i)(k-i)}}_{+}^{\prime}(x)\bigg)
\end{align*}
for each $x\in (a,b)$. By our assumption on $f$ we have $f^{\prime}\geq 0$ and $f^{(n-k-1)}\geq 0$. We also have $\paren*{\prod_{i=1}^{k}h_i^{\varphi(i)(k-i)}}_{+}^{\prime}\geq 0$ 
because of (iii) and 
$$\paren*{\paren*{\prod_{i\in I}h_i^i}\paren*{\prod_{i=k+1}^nh_i^{\varphi(i)(i-k)}}}_{+}^{\prime}\leq 0$$ 
because of (ii).\par
Furthermore, the assumption $\alpha\leq\frac{n}{n-1}$ ensures that $n+1>\alpha$ and $n(1-\alpha)+k\alpha\geq 0$. It follows that, for all $x\in (a,b)$,
\begin{align*}
&(H_{k+1})_{+}^{\prime}(x)\leq \\
&f^{\alpha}(x)\paren*{\prod_{i=k}^ni}\paren*{\int_a^xf^{\alpha}\prod_{j\in J}h_j}^{k-1}\paren*{\prod_{j\in J}h_j(x)}\paren*{\prod_{i\in I}h_i^i(x)}
\paren*{\prod_{i=k+1}^nh_i^{\varphi(i)(i-k)}(x)} \\
&-(n+1-\alpha)f^{n(1-\alpha)+k\alpha}(x)f^{(n-k)}(x)\prod_{i=1}^{k}h_i^{\varphi(i)(k-i)}(x) \\
&=f^{\alpha}(x)H_k(x)\prod_{i=1}^{k-1}h_i^{\varphi(i)}(x)\leq 0.
\end{align*}
Thus $H_{k+1}$ is decreasing and hence $H_{k+1}(x)\leq H_{k+1}(a)\leq 0$ for all $x\in [a,b]$, which finishes the induction.\par
Now we define
\begin{equation*}
H(x):=nh_n^{\varphi(n)}(x)\paren*{\int_a^xf^{\alpha}\prod_{j\in J}h_j}^{n-1}\paren*{\prod_{i\in I}h_i^i(x)}-f^{n+1-\alpha}(x)\prod_{i=1}^{n-1}h_i^{\varphi(i)(n-i-1)}(x).
\end{equation*}
Then $H\in D_{+}[a,b]$ and similar to the induction step above one can show that
\begin{equation*}
H_{+}^{\prime}(x)\leq f^{\alpha}(x)H_{n-1}(x)\prod_{i=1}^{n-2}h_i^{\varphi(i)}(x) \ \ \ \forall x\in (a,b).
\end{equation*}
Since $H_{n-1}\leq 0$ it follows that $H$ is decreasing and hence $H(x)\leq H(a)\leq 0$ for every $x\in [a,b]$. Also, from the definition $G$ one can easily see 
that $G=H\prod_{i=1}^{n-1}h_i^{\varphi(i)}$. Together with \eqref{Eq:FG} we obtain
\begin{equation*}
F^{\prime}(x)\leq f^{\alpha}(x)H(a)M^n\paren*{\prod_{i\in I}m_i^{n-i}}\paren*{\prod_{i=1}^{n-1}m_i^{\varphi(i)}} \ \ \forall x\in [a,b]
\end{equation*}
Furthermore, it is easily checked that 
\begin{equation*}
H(a)\prod_{i=1}^{n-1}m_i^{\varphi(i)}\leq -f^{n+1-\alpha}(a)\prod_{i\in J}m_i^{n-i}.
\end{equation*}
Also, $f$ is increasing (since $f^{\prime}\geq 0$) and hence $f^{\alpha}$ is increasing for $\alpha\geq 0$ and decreasing for $\alpha<0$.
Thus we get
\begin{equation*}
F^{\prime}(x)\leq -KM^n\prod_{i=1}^nm_i^{n-i} \ \ \ \forall x\in [a,b].
\end{equation*}
By the mean value theorem this implies
\begin{equation*}
F(b)=F(b)-F(a)\leq -(b-a)KM^n\prod_{i=1}^nm_i^{n-i},
\end{equation*}
which is the desired inequality.
\end{Proof}

Let us explicitly note the following special case of Proposition \ref{prop:general-qi-inequal}.
\begin{corollary}\label{cor:fghp}
Let $g,h,p$ be nonnegative functions on $[a,b]$ such that $h, p\in D_{+}[a,b]$ 
or $h, p\in D_{-}[a,b]$, $g$ is continuous, $p$ is decreasing and $h$ and 
$gp$ are increasing.\par
Suppose further that $n\in \mathbb{N}$, $n\geq 2$, and $f$ is a strictly positive, $(n-1)$-times differentiable function on $[a,b]$ 
with $f^{(i)}\geq 0$ for $i=1,\dots,n-1$ and $(n+1-\alpha)f^{(n-1)}f^{n(1-\alpha)}\geq n!hp^{\nu}$ for some $\nu\in \set*{2,\dots,n}$ and $\alpha\leq \frac{n}{n-1}$. 
Let $M:=\inf\set*{g(x):x\in [a,b]}$. Then we have
\begin{equation*}
(b-a)KM^nh^{n-1}(a)p^{n-\nu}(b)+\paren*{\int_a^bf^{\alpha}ghp}^n\leq\int_a^bf^{n+1}g^nh^{n-1}p^{n-\nu},
\end{equation*}
where $K=f^{n+1}(a)$ if $\alpha\geq 0$ and $K=f^{n+1-\alpha}(a)f^{\alpha}(b)$ if $\alpha<0$.
\end{corollary}

\begin{Proof}
Put $h_1:=h$, $h_{\nu}:=p$ and $h_i:=1$ for $i\in \set*{2,\dots,n}\sm \set*{\nu}$ as well as $I:=\set*{\nu}, J:=\set*{1,\dots,n}\sm I$ 
and apply Proposition \ref{prop:general-qi-inequal}.
\end{Proof}

This yields in particular the following Corollaries.

\begin{corollary}\label{cor:fg1}
Let $\alpha\leq \frac{n}{n-1}$. If $f$ is a strictly positive, $(n-1)$-times differentiable function on $[a,b]$ (where $n\geq 2$) such that 
$f^{(i)}\geq 0$ for $i=1,\dots,n-1$, then 
\begin{equation*}
(b-a)Kg^n(a)+(n+1-\alpha)\frac{A}{n!}\paren*{\int_a^bf^{\alpha}g}^n\leq \int_a^bf^{n+1}g^n
\end{equation*}
holds for every nonnegative function $g$ on $[a,b]$ which is continuous and increasing, where $A:=\inf\set*{f^{(n-1)}(x)f^{n(1-\alpha)}(x):x\in [a,b]}$
and $K$ is defined as before.
\end{corollary}

\begin{Proof}
Put $c:=((n+1-\alpha)A/n!)^{-1/(n(1-\alpha)+1)}$ and apply Corollary \ref{cor:fghp} to the function $cf$ (with $h=p=1$).
\end{Proof}

\begin{corollary}\label{cor:fg2}
Let $f$ be a strictly positive function on $[a,b]$ which is $n$-times differentiable $(n\geq 2)$ with $f^{(i)}\geq 0$ for all $i=1,\dots,n$. 
Let $g$ be a nonnegative function on $[a,b]$ which is continuous and increasing and let $\alpha\leq \frac{n}{n-1}$. Then we have
\begin{align*}
&(b-a)Kg^n(a)(f^{(n-1)}(a))^{n-1}+(n+1-\alpha)\frac{C^{n(1-\alpha)}}{n!}\paren*{\int_a^bf^{\alpha}gf^{(n-1)}}^n\\
&\leq\int_a^bf^{n+1}g^n(f^{(n-1)})^{n-1},
\end{align*}
where $K$ is defined as above and $C=f(a)$ if $\alpha\leq 1$, $C=f(b)$ if $\alpha>1$.
\end{corollary}

\begin{Proof}
Set $p:=1$, $h:=(n+1-\alpha)C^{n(1-\alpha)}f^{(n-1)}/n!$ and $\nu:=2$ and apply Corollary \ref{cor:fghp}.
\end{Proof}

\begin{corollary}\label{cor:fh}
Let $f$ be a strictly positive, $n$-times differentiable function on $[a,b]$ $(n\geq 2)$ satisfying $f^{(i)}\geq 0$ 
for $i=1,\dots,n-2$, $f^{(n-1)}(x)>0$ for all $x\in [a,b]$, and $f^{(n)}\leq 0$. Let $h\in D_{+}[a,b]$ or $h\in D_{-}[a,b]$ 
be increasing with $0\leq h\leq 1$ and let $\alpha\leq \frac{n}{n-1}$. Then we have
\begin{equation*}
(b-a)K\frac{h^{n-1}(a)}{f^{(n-1)}(a)}+(n+1-\alpha)\frac{C^{n(1-\alpha)}}{n!}\paren*{\int_a^bf^{\alpha}h}^n\leq \int_a^b\frac{f^{n+1}h^{n-1}}{f^{(n-1)}},
\end{equation*}
where $K$ and $C$ are defined as in the previous Corollary.
\end{corollary}

\begin{Proof}
Put $p(t):=((n+1-\alpha)C^{n(1-\alpha)}f^{(n-1)}(t)/n!)^{1/n}$ for $t\in [a,b]$, $g:=1/p$ and $\nu:=n$ and apply Corollary \ref{cor:fghp}.
\end{Proof}

\begin{corollary}\label{cor:fg3}
Let $f$ be a strictly positive, $n$-times differentiable function on $[a,b]$ $(n\geq 2)$ satisfying $f^{(i)}\geq 0$ for $i=1,\dots,n$.
Let $g$ be a continuous, increasing function on $[a,b]$ with $g\geq 0$ and let $\alpha\in (-\infty,1]$. Put $\beta:=n(1-\alpha)+\alpha$.
Then we have
\begin{align*}
&(b-a)Lg^n(a)(f^{(n-1)}(a))^{n-1}+\frac{n+1-\alpha}{n!}\paren*{\int_a^bf^{\beta}gf^{(n-1)}}^n \\
&\leq\int_a^bf^{n\beta+1}g^n(f^{(n-1)})^{n-1},
\end{align*}
where $L=f^{n\beta+1}(a)$ for $\alpha\geq 0$ and $L=f^{n\beta+1-\alpha}(a)f^{\alpha}(b)$ for $\alpha<0$.
\end{corollary}

\begin{Proof}
Put $h:=(n+1-\alpha)f^{(n-1)}f^{n(1-\alpha)}/n!$, $p:=1$ and $\nu:=2$ and apply Corollary \ref{cor:fghp}.
\end{Proof}

\section{An integral inequality for 1/f}\label{sec:g-over-f}

Proposition \ref{prop:general-qi-inequal} already provides a Qi-type integral inequality for $1/f$ (consider the case $\alpha=-1$).
The next result is another inequality in the spirit of \eqref{eq:qi} for $1/f$.
\begin{proposition}\label{prop:g-over-f}
Let $f:[a,b] \rightarrow \R$ be twice differentiable such that $f^{\prime}, f^{\prime\prime}\geq 0$ and $f(t)>0$ for all $t\in [a,b]$. Let further 
$g$ be an increasing, continuous function on $[a,b]$ with $g\geq 0$ and $h\in D_{+}[a,b]$ or $h\in D_{-}[a,b]$ be an increasing, nonnegative function.\par 
If $n\in \mathbb{N}$ and $f^{n+1}(a)(f^{\prime}(a))^nh(a)\geq n!/(n+1)^{n-1}$, then
\begin{equation*}
(b-a)\frac{f^{n+1}(a)g^{n+1}(a)h(a)}{f(b)}+\paren*{\int_a^b\frac{g}{f}}^{n+1}\leq\int_a^bf^ng^{n+1}h.
\end{equation*}
\end{proposition}

Note that the conditions $f^{\prime}\geq 0$ and $f^{\prime\prime}\geq 0$ just mean that $f$ is increasing and convex.

\begin{Proof}
Let $h\in D_{+}[a,b]$ (the other case is completely analogous).
We define
\begin{equation*}
F(x):=\paren*{\int_a^x\frac{g}{f}}^{n+1}-\int_a^xf^ng^{n+1}h \ \ \forall x\in [a,b]
\end{equation*}
and
\begin{equation*}
G(x):=(n+1)\paren*{\int_a^x\frac{1}{f}}^n-f^{n+1}(x)h(x) \ \ \forall x\in [a,b]
\end{equation*}
and claim that
\begin{equation}\label{eq:2:FG}
F^{\prime}\leq \frac{g^{n+1}}{f}G.
\end{equation}
To see this first note that
\begin{equation}\label{eq:2:diffF}
F^{\prime}(x)=\frac{g(x)}{f(x)}\paren*{(n+1)\paren*{\int_a^x\frac{g}{f}}^n-f^{n+1}(x)g^n(x)h(x)}.
\end{equation}
Since $g$ is increasing we have
\begin{equation*}
\int_a^x\frac{g}{f}\leq g(x)\int_a^x\frac{1}{f}.
\end{equation*}
Combining this with \eqref{eq:2:diffF} gives \eqref{eq:2:FG}.\par
We further have $G\in D_{+}[a,b]$ and 
\begin{equation*}
G_{+}^{\prime}(x)=\frac{n(n+1)}{f(x)}\paren*{\int_a^x\frac{1}{f}}^{n-1}-(n+1)f^n(x)f^{\prime}(x)h(x)-f^{n+1}(x)h_{+}^{\prime}(x).
\end{equation*}
Since $h$ is increasing we have $h_{+}^{\prime}\geq 0$ and hence 
\begin{equation}\label{eq:2:diffG}
G_{+}^{\prime}(x)\leq\frac{n+1}{f(x)}\paren*{n\paren*{\int_a^x\frac{1}{f}}^{n-1}-f^{n+1}(x)f^{\prime}(x)h(x)}.
\end{equation}
Next we define functions $H_1,\dots,H_n$ on $[a,b]$ by
\begin{equation*}
H_k(x):=\paren*{\prod_{i=k}^ni}\paren*{\int_a^x\frac{1}{f}}^{k-1}-(n+1)^{n-k}f^{n+1}(x)(f^{\prime})^{n-k+1}(x)h(x).
\end{equation*}
We will show inductively that $H_k\leq 0$ for all $k=1,\dots,n$. First note that
$H_1=n!-(n+1)^{n-1}f^{n+1}(f^{\prime})^nh$ and by our assumptions $h$, $f$ and $f^{\prime}$ are increasing functions, thus $H_1$ is decreasing.
Hence $H_1\leq H_1(a)$ and again by assumption we have $H_1(a)\leq 0$.\par
Now suppose that $1\leq k<n$ and $H_k\leq 0$. We have $H_{k+1}\in D_{+}[a,b]$ and 
\begin{align*}
&(H_{k+1})_{+}^{\prime}(x)=\frac{1}{f(x)}\paren*{\prod_{i=k}^ni}\paren*{\int_a^x\frac{1}{f}}^{k-1}-(n+1)^{n-k}f^n(x)(f^{\prime})^{n-k+1}(x)h(x) \\
&-(n+1)^{n-k-1}f^{n+1}(x)((f^{\prime})^{n-k}h)_{+}^{\prime}(x).
\end{align*}
Since $f^{\prime}$ and $h$ are increasing, it follows that
\begin{equation*}
(H_{k+1})_{+}^{\prime}(x)\leq\frac{H_k(x)}{f(x)}\leq 0 \ \ \forall x\in (a,b)
\end{equation*}
and hence $H_{k+1}\leq H_{k+1}(a)\leq 0$.\par
So in particular $H_n\leq 0$ and from \eqref{eq:2:diffG} it follows that $G_{+}^{\prime}\leq (n+1)H_n/f$. Thus $G_{+}^{\prime}\leq 0$ and consequently,
$G\leq G(a)=-f^{n+1}(a)h(a)$.\par
Using \eqref{eq:2:FG} and the fact that $g$ is increasing and $1/f$ decreasing, we obtain $F^{\prime}\leq -f^{n+1}(a)g^{n+1}(a)h(a)/f(b)$. 
The mean value theorem now implies $F(b)=F(b)-F(a)\leq -(b-a)f^{n+1}(a)g^{n+1}(a)h(a)/f(b)$, which is equivalent to the desired inequality.
\end{Proof}

Let us now collect some corollaries to the above result.
\begin{corollary}\label{cor:g-over-f}
Let $f$ be a twice differentiable, strictly positive function on $[a,b]$ such that $f^{\prime}(a)>0$ and $f^{\prime\prime}\geq 0$ and
let $g$ be a nonnegative function on $[a,b]$ which is continuous and increasing. Then we have for each $n\in \N$
\begin{align*}
(b-a)\frac{g^{n+1}(a)}{f(b)}+(f^{\prime}(a))^n\frac{(n+1)^{n-1}}{n!}\paren*{\int_a^b\frac{g}{f}}^{n+1}\leq\frac{1}{f^{n+1}(a)}\int_a^bf^ng^{n+1}.
\end{align*}
\end{corollary}

\begin{Proof}
For a given $n\in \N$, put $c:=(n!((n+1)^{n-1}f^{n+1}(a)(f^{\prime}(a))^n)^{-1})^{1/(2n+1)}$ and apply Proposition \ref{prop:g-over-f} to the 
functions $cf$ and $g$ (and $h:=1$).
\end{Proof}

\begin{corollary}\label{cor:limit}
If $f$ and $g$ are as in the previous corollary, then 
\begin{align*}
\int_a^b\frac{g}{f}&\leq\frac{1}{ef(a)f^{\prime}(a)}\liminf_{n\to \infty}\paren*{\int_a^bf^ng^{n+1}}^{1/(n+1)} \\
&\leq\frac{1}{ef(a)f^{\prime}(a)}\limsup_{n\to \infty}\paren*{\int_a^bf^ng^{n+1}}^{1/(n+1)}\leq\frac{f(b)g(b)}{ef(a)f^{\prime}(a)}.
\end{align*}
\end{corollary}

\begin{Proof}
Since $f$ and $g$ are increasing, we have
\begin{equation*}
\paren*{\int_a^bf^ng^{n+1}}^{1/(n+1)}\leq (b-a)^{1/(n+1)}(f(b))^{n/(n+1)}g(b) \ \ \forall n\in \N.
\end{equation*}
The righthand side of this inequality converges to $f(b)g(b)$ for $n\to \infty$, thus
\begin{equation*}
\limsup_{n\to \infty}\paren*{\int_a^bf^ng^{n+1}}^{1/(n+1)}\leq f(b)g(b).
\end{equation*}
Corollary \ref{cor:g-over-f} further implies that
\begin{equation*}
\int_a^b\frac{g}{f}\leq\paren*{\int_a^bf^ng^{n+1}}^{1/(n+1)}\paren*{\frac{n!}{(n+1)^{n-1}}}^{1/(n+1)}\frac{1}{f(a)(f^{\prime}(a))^{n/(n+1)}}
\end{equation*}
for each $n\in \N$.\par
Using the well-known limits $\lim_{n\to \infty}\sqrt[n]{n}=1$ and $\lim_{n\to \infty}n/\sqrt[n]{n!}=e$, we obtain
\begin{equation*}
\lim_{n\to \infty}\paren*{\frac{n!}{(n+1)^{n-1}}}^{1/(n+1)}=\lim_{n\to \infty}\frac{\sqrt[n+1]{(n+1)!}}{n+1}\sqrt[n+1]{n+1}=\frac{1}{e}.
\end{equation*}
It follows that
\begin{equation*}
\int_a^b\frac{g}{f}\leq\frac{1}{ef(a)f^{\prime}(a)}\liminf_{n\to \infty}\paren*{\int_a^bf^ng^{n+1}}^{1/(n+1)}.
\end{equation*}
\end{Proof}

In fact, the inequality $\int_a^b\frac{g}{f}\leq f(b)g(b)/(ef(a)f^{\prime}(a))$
also holds under weaker assumptions.
\begin{lemma}\label{lemma:1/e}
Let $f:[a,b] \rightarrow \R$ be a strictly positive, differentiable function such
that $f^{\prime}$ is increasing and $f^{\prime}(a)>0$. Let $g:[a,b] \rightarrow \R$
be increasing and nonnegative. Then we have
\begin{equation}
\int_a^b\frac{g}{f}\leq \frac{g(b)}{f^{\prime}(a)}\log\paren*{\frac{f(b)}{f(a)}}
\leq\frac{f(b)g(b)}{ef(a)f^{\prime}(a)}.
\end{equation}
\end{lemma}

For the proof we need the following Lemma (which is surely well-known, but the author was unable to find a reference).
\begin{lemma}\label{lemma:exp}
For every $x>0$ we have $x^e\leq e^x$. Equality holds if and only if $x=e$.
\end{lemma}

\begin{Proof}
Put $h(x):=x-e\log(x)$ for $x>0$. Then $h^{\prime}(x)=1-e/x$ and hence
$h^{\prime}(x)>0$ for $x>e$ and $h^{\prime}(x)<0$ for $x<e$.\par 
Thus $h$ is strictly increasing on $[e,\infty)$ and strictly decreaisng on $(0,e]$.
This implies $h(x)>h(e)=0$ for all $x>0$ with $x\neq e$, which implies the claimed
inequality.
\end{Proof}

\begin{Proof}{(of Lemma \ref{lemma:1/e})}
The monotonicity of $g$ and $f^{\prime}$ implies
\begin{equation*}
\int_a^b\frac{g}{f}\leq \frac{g(b)}{f^{\prime}(a)}\int_a^b\frac{f^{\prime}}{f}
=\frac{g(b)}{f^{\prime}(a)}\log\paren*{\frac{f(b)}{f(a)}}
\end{equation*}	
and Lemma \ref{lemma:exp} implies 
\begin{equation*}
\log\paren*{\frac{f(b)}{f(a)}}\leq \frac{f(b)}{ef(a)}
\end{equation*}
which concludes the proof.
\end{Proof}

Next we will derive three more corollaries concerning the logarithm of a function $f$.
\begin{corollary}\label{cor:log-1}
Let $f$ be as in Corollary \ref{cor:g-over-f} and $n\in \N$. Then we have
\begin{equation*}
\paren*{\log\paren*{\frac{f(b)}{f(a)}}}^{n+1}+\frac{n!(b-a)}{(n+1)^{n-1}}\frac{f^{\prime}(a)}{f(b)}\leq
\frac{n!}{(n+1)^n}\frac{(f^{\prime}(b))^n}{(f^{\prime}(a))^n}\paren*{\frac{f^{n+1}(b)}{f^{n+1}(a)}-1}.
\end{equation*}
\end{corollary}

\begin{Proof}
We apply Corollary \ref{cor:g-over-f} with $g:=f^{\prime}$ to get
\begin{align*}
&\paren*{\log\paren*{\frac{f(b)}{f(a)}}}^{n+1}=\paren*{\int_a^b\frac{f^{\prime}}{f}}^{n+1} \\
&\leq\frac{n!}{(n+1)^{n-1}f^{n+1}(a)(f^{\prime}(a))^n}\int_a^bf^n(f^{\prime})^{n+1}-(b-a)\frac{n!f^{\prime}(a)}{f(b)(n+1)^{n-1}}.
\end{align*}
Since $f^{\prime}\leq f^{\prime}(b)$ we have
\begin{equation*}
\int_a^bf^n(f^{\prime})^{n+1}\leq (f^{\prime}(b))^n\int_a^bf^nf^{\prime}=\frac{(f^{\prime}(b))^n}{n+1}(f^{n+1}(b)-f^{n+1}(a)).
\end{equation*}
Combining these two estimates and simplifying a little finishes the proof.
\end{Proof}

\begin{corollary}\label{cor:log-2}
Let $f$ be as in Corollary \ref{cor:g-over-f} and assume in addition that $\log(f)$ is convex. Then we have for each $n\in \N$
\begin{equation*}
\paren*{1-\frac{f(a)}{f(b)}}^{n+1}\leq\frac{n!}{(n+1)^{n-1}}\paren*{\frac{(f^{\prime}(b))^n}{(f^{\prime}(a))^n}\log\paren*{\frac{f(b)}{f(a)}}-(b-a)\frac{f^{\prime}(a)}{f(b)}}
\end{equation*}
\end{corollary}

\begin{Proof}
Since $\log(f)$ is convex, the derivative $(\log(f))^{\prime}=f^{\prime}/f$ is increasing. Thus we can apply Corollary \ref{cor:g-over-f} with $g:=f^{\prime}/f$ to get
\begin{align*}
&\paren*{\frac{1}{f(a)}-\frac{1}{f(b)}}^{n+1}=\paren*{\int_a^b\frac{f^{\prime}}{f^2}}^{n+1} \\
&\leq \frac{n!}{(n+1)^{n-1}}\paren*{\frac{1}{f^{n+1}(a)(f^{\prime}(a))^n}\int_a^b\frac{(f^{\prime})^{n+1}}{f}-(b-a)\frac{f^{\prime}(a)}{f(b)f^{n+1}(a)}}.
\end{align*}
Since $f^{\prime}$ is increasing we can estimate
\begin{equation*}
\int_a^b\frac{(f^{\prime})^{n+1}}{f}\leq (f^{\prime}(b))^n\int_a^b\frac{f^{\prime}}{f}=(f^{\prime}(b))^n\log\paren*{\frac{f(b)}{f(a)}}.
\end{equation*}
Combining the two estimates and multiplying by $(f(a))^{n+1}$ gives the desired inequality.
\end{Proof}

Applying Lemma \ref{lemma:1/e} to $g:=f^{\prime}/f$ also yields the following result.
\begin{corollary}\label{cor:log-3}
If $f\in C^1[a,b]$ is a strictly positive function with a strictly positive derivative
such that $\log(f)$ is convex, then 
\begin{equation*}
1-\frac{f(a)}{f(b)}\leq \frac{f^{\prime}(b)f(a)}{f(b)f^{\prime}(a)}\log\paren*{\frac{f(b)}{f(a)}}\leq \frac{f^{\prime}(b)}{ef^{\prime}(a)}.
\end{equation*}
\end{corollary} 

\section{Further variants of Qi's inequality}\label{sec:alpha}

In the last section of this paper, we present some other variants of Qi's inequality \eqref{eq:qi}. We start with some results in which we have the same 
exponent $n$ for the integral and the function $f$.

\begin{proposition}\label{prop:alpha1}
Let $f$ be an $(n-1)$-times differentiable function on $[a,b]$ ($n\geq 2$) such that $f^{(i)}\geq 0$ for $i=0,\dots,n-2$.
Let further $g\in D_{+}[a,b]$ or $g\in D_{-}[a,b]$ be a strictly positive, increasing function and let $\alpha\in (n,\infty)$.
Suppose that $g^{\alpha-n}f^{(n-1)}\geq \frac{n!}{n-1}f$. Then we have
\begin{equation*}
(b-a)f^n(a)g^{\alpha}(a)+\paren*{\int_a^bfg}^n\leq\int_a^bf^ng^{\alpha}.
\end{equation*}
\end{proposition}

\begin{Proof}
The proof is similar to the previous ones. Again we assume that $g\in D_{+}[a,b]$, define
\begin{align*}
&F(x):=\paren*{\int_a^xfg}^n-\int_a^xf^ng^{\alpha}, \\
&G(x):=n\paren*{\int_a^xf}^{n-1}-f^{n-1}(x)g^{\alpha-n}(x)
\end{align*}
and show as before that
\begin{equation}\label{eq:Fprime}
F^{\prime}\leq fg^nG.
\end{equation}
Next, using the fact that $g_{+}^{\prime}\geq 0$ and $\alpha>n$, we obtain
\begin{equation}\label{eq:Gprime}
G_{+}^{\prime}(x)\leq (n-1)\paren*{nf(x)\paren*{\int_a^xf}^{n-2}-f^{n-2}(x)f^{\prime}(x)g^{\alpha-n}(x)}.
\end{equation}
If $n=2$ this implies $G_{+}^{\prime}\leq 2f-f^{\prime}g^{\alpha-2}$, which by assumption is nonpositive. Thus we obtain $G\leq G(a)=-f(a)g^{\alpha-2}(a)$.\par
Now consider the case $n\geq 3$. We define functions $H_1,\dots, H_{n-2}$ on $[a,b]$ by
\begin{equation*}
H_k(x):=n\paren*{\prod_{i=k+1}^{n-2}i}\paren*{\int_a^xf}^{k}-f^{k-1}(x)f^{(n-k-1)}(x)g^{\alpha-n}(x).
\end{equation*}
We will show that $H_k\leq 0$ for $k=1,\dots,n-2$. First, we have $(H_1)_{+}^{\prime}=fn!/(n-1)-f^{(n-1)}g^{\alpha-n}-f^{(n-2)}(\alpha-n)g^{\alpha-n-1}g_{+}^{\prime}$.
Since $f^{(n-2)}, g_{+}^{\prime}\geq 0$ and $\alpha>n$ this implies $(H_1)_{+}^{\prime}\leq fn!/(n-1)-f^{(n-1)}g^{\alpha-n}$, 
and thus by our assumption on $f$ and $g$ we have $(H_1)_{+}^{\prime}\leq 0$. Hence $H_1\leq H_1(a)\leq 0$.\par
Now suppose that $H_k\leq 0$ for some $k\leq n-3$. We have
\begin{align*}
&(H_{k+1})_{+}^{\prime}(x)=nf(x)\paren*{\prod_{i=k+1}^{n-2}i}\paren*{\int_a^xf}^{k}-kf^{k-1}(x)f^{\prime}(x)f^{(n-k-2)}(x)g^{\alpha-n}(x) \\
&-f^{k}(x)f^{(n-k-1)}(x)g^{\alpha-n}(x)-f^{k}(x)f^{(n-k-2)}(x)(\alpha-n)g^{\alpha-n-1}(x)g_{+}^{\prime}(x) \\
&\leq nf(x)\paren*{\prod_{i=k+1}^{n-2}i}\paren*{\int_a^xf}^{k}-f^{k}(x)f^{(n-k-1)}(x)g^{\alpha-n}(x)=f(x)H_k(x)\leq 0.
\end{align*}
Hence $H_{k+1}\leq H_{k+1}(a)\leq 0$.\par
It follows from \eqref{eq:Gprime} that $G_{+}^{\prime}\leq (n-1)fH_{n-2}\leq 0$ and thus $G\leq G(a)=-f^{n-1}(a)g^{\alpha-n}(a)$.\par
Using this together with \eqref{eq:Fprime} and the fact that $f$ and $g$ are increasing we obtain $F^{\prime}\leq -f^n(a)g^{\alpha}(a)$.
The mean value theorem therefore implies $F(b)=F(b)-F(a)\leq -(b-a)f^n(a)g^{\alpha}(a)$, finishing the proof.
\end{Proof}

This yields the following corollary.

\begin{corollary}\label{cor:alpha1}
If $n\geq 2$ and $f$ is an $(n-1)$-times differentiable function on $[a,b]$ such that $f^{(i)}\geq 0$ for $i=0,\dots,n-2$, $f^{(n-1)}$ is strictly positive
and $f/f^{(n-1)}$ is bounded, then
\begin{equation*}
(b-a)f^n(a)+\frac{n-1}{n!A}\paren*{\int_a^bf}^n\leq \int_a^bf^n,
\end{equation*}
where $A:=\norm*{f/f^{(n-1)}}_{\infty}$ and $\norm{\cdot}_{\infty}$ denotes the usual sup-norm.
\end{corollary}

\begin{Proof}
Define $g$ to be constant on $[a,b]$ with value $\paren*{n!A/(n-1)}^{1/n}$ and $\alpha:=2n$ and apply
Proposition \ref{prop:alpha1}.
\end{Proof}

An analogous result to Proposition \ref{prop:alpha1} for the case $\alpha\leq n$ reads as follows.

\begin{proposition}\label{prop:alpha2}
Let $f$ be an $(n-1)$-times differentiable function on $[a,b]$ ($n\geq 2$) such that $f^{(i)}\geq 0$ for $i=0,\dots,n-2$.
Let $g$ be a strictly positive, continuous and increasing function on $[a,b]$ and let $\alpha\in (-\infty,n]$ with
$f^{(n-1)}\geq \frac{n!}{n-1}f\norm{g}_{\infty}^{n-\alpha}$. Then we have
\begin{equation*}
(b-a)f^n(a)K+\paren*{\int_a^bfg}^n\leq\int_a^bf^ng^{\alpha},
\end{equation*}
where $K=g^{\alpha}(a)$ for $\alpha\geq 0$ and $K=g^{\alpha}(b)$ for $\alpha<0$.
\end{proposition}

\begin{Proof}
The proof is similar to the last one, so we will only sketch it. First define $F$ exactly as in the previous proof
and let
\begin{equation*}
G(x):=n\norm{g}_{\infty}^{n-\alpha}\paren*{\int_a^x f}^{n-1}-f^{n-1}(x).
\end{equation*}
It is easy to prove that
\begin{equation}\label{eq:FprimeG}
F^{\prime}\leq fg^{\alpha}G.
\end{equation}
Next we define again functions $H_k$ (where $k=1,\dots,n-2$) by
\begin{equation*}
H_k(x):=n\norm{g}_{\infty}^{n-\alpha}\paren*{\prod_{i=k+1}^{n-2}i}\paren*{\int_a^xf}^k-f^{k-1}(x)f^{(n-k-1)}(x).
\end{equation*}
Similar to our previous proofs one can show inductively that $H_k\leq 0$ for $k=1,\dots,n-2$ and $G^{\prime}=(n-1)fH_{n-2}$.\par
It follows that $G$ is decreasing and hence $G\leq G(a)=-f^{n-1}(a)$. Using this together with
\eqref{eq:FprimeG}, the monotonicity of $f$ and $g$, and the mean value theorem, one obtains
\begin{equation*}
F(b)\leq -(b-a)f^n(a)K
\end{equation*}
and the proof is finished.
\end{Proof}

Here is yet another result of the above type.
\begin{proposition}\label{prop:fgh}
Let $f$ be an $(n-1)$-times differentiable function on $[a,b]$ ($n\geq 2$) such that $f^{(i)}\geq 0$ for $i=0,\dots,n-2$.
Let $g, h\in D_{+}[a,b]$ or $g, h\in D_{-}[a,b]$ with $g, h\geq 0$. Assume further that $g$ is increasing and $h$ and $gh$ are decreasing. 
Let $f^{(n-1)}\geq\frac{n!}{n-1}fg^{n-l}h^n$ for some $l\in \set*{1,\dots,n}$.
Then we have
\begin{equation*}
(b-a)f^n(a)g^l(a)+\paren*{\int_a^bfgh}^n\leq \int_a^bf^ng^l.
\end{equation*}
\end{proposition}

\begin{Proof}
We will only sketch the proof for the case $g, h\in D_{+}[a,b]$. First define
\begin{align*}
&F(x):=\paren*{\int_a^xfgh}^n-\int_a^xf^ng^l, \\
&G(x):=n\paren*{\int_a^xfgh}^{n-1}h(x)-f^{n-1}(x)g^{l-1}(x).
\end{align*}
Then $F^{\prime}=fgG$.\par
In the following, we will only treat the case $2\leq l\leq n-1$. The boundary cases are treated similarly.
We define
\begin{align*}
&H_k(x):=n\paren*{\prod_{i=2}^k(n-i)}\paren*{\int_a^xfgh}^{n-k-1}h^{k+1}(x)-f^{(k)}(x)f^{n-k-2}(x)g^{l-k-1}(x), \\
&J_s(x):=n\paren*{\prod_{i=2}^{l+s-1}(n-i)}\paren*{\int_a^xfgh}^{n-l-s}h^{l+s}(x)g^s(x)-f^{(l+s-1)}(x)f^{n-l-s-1}(x)
\end{align*}
for $k=1,\dots,l-1$ and $s=0,\dots,n-1-l$.\par 
Then $J_0=H_{l-1}$. Moreover, using $h^{\prime}\leq 0$ and $g^{\prime}\geq 0$, it is easy to see that $G_{+}^{\prime}\leq (n-1)fgH_1$,
$(H_k)_{+}^{\prime}\leq fgH_{k+1}$ for each $k$ and $(J_s)_{+}^{\prime}\leq fJ_{s+1}$ for each $s$ (note that $h^{l+s}g^s$ is also decreasing
(and hence has a negative right-derivative) since $h^{l+s}g^s=(gh)^sh^l$).\par
We further have $(J_{n-l-1})_{+}^{\prime}\leq \frac{n!}{n-1}fh^ng^{n-l}-f^{(n-1)}$, whence by assumption $(J_{n-l-1})_{+}^{\prime}\leq 0$.\par 
Using the usual monotonicity and induction arguments, we obtain that all the functions $J_s$ and $H_k$ are negative. It follows
in particular that $G_{+}^{\prime}\leq 0$ and hence $G\leq G(a)=-f^{n-1}(a)g^{l-1}(a)$.\par
Since $F^{\prime}=fgG$ and $f$ and $g$ are increasing, we obtain $F^{\prime}\leq -f^n(a)g^l(a)$. From this and the mean value theorem
we can deduce the desired inequality.
\end{Proof}

An immediate corollary is the following inequality.
\begin{corollary}\label{cor:fgh1}
Let $f$ be an $(n-1)$-times differentiable function on $[a,b]$ ($n\geq 2$) such that $f^{(i)}\geq 0$ for $i=0,\dots,n-2$ 
and $g$ a strictly positive, increasing function on $[a,b]$ with $g\in D_{+}[a,b]$ or $g\in D_{-}[a,b]$. If $f^{(n-1)}g\geq\frac{n!}{n-1}f$, then
\begin{equation*}
(b-a)f^n(a)g(a)+\paren*{\int_a^bf}^n\leq \int_a^bf^ng.
\end{equation*}
\end{corollary}

\begin{Proof}
Put $h:=1/g$ and $l=1$ and apply Proposition \ref{prop:fgh}.
\end{Proof}

This yields a further corollary which is related to Corollary \ref{cor:fh}.
\begin{corollary}\label{cor:fgh2}
Let $f$ be a strictly positive, $n$-times differentiable function on $[a,b]$ ($n\geq 2$) such that $f^{(i)}\geq 0$ for $i=1,\dots,n-2$.
Suppose further that $f^{(n-1)}$ is strictly positive and $f/f^{(n-1)}$ is increasing. Then we have
\begin{equation*}
(b-a)\frac{f^{n+1}(a)}{f^{(n-1)}(a)}+\frac{n-1}{n!}\paren*{\int_a^bf}^n\leq \int_a^b\frac{f^{n+1}}{f^{(n-1)}}.
\end{equation*}
\end{corollary}

\begin{Proof}
Put $g:=n!/(n-1)(f/f^{(n-1)})$ and apply Corollary \ref{cor:fgh1}.
\end{Proof}

Note that Corollary \ref{cor:fh} implies (for $h=1$ and $\alpha=1$) the stronger inequality
\begin{equation*}
(b-a)\frac{f^{n+1}(a)}{f^{(n-1)}(a)}+\frac{1}{(n-1)!}\paren*{\int_a^bf}^n\leq \int_a^b\frac{f^{n+1}}{f^{(n-1)}}
\end{equation*}
under the stronger assumption that $f^{(n)}\leq 0$ (the assumption that $f/f^{(n-1)}$ is increasing is equivalent to $f^{(n)}\leq f^{(n-1)}f^{\prime}/f$).

Finally, we also have the following generalisation of \eqref{eq:qi}.

\begin{proposition}\label{prop:x-minus-a}
Let $n\geq 2$, $k\in \set*{1,\dots,n-1}$ and $f:[a,b] \rightarrow \R$ a $k$-times differentiable function
such that $f^{(i)}\geq 0$ for $i=0,\dots,k-1$. Let $g:[a,b] \rightarrow \R$ be continuous and increasing with $g\geq 0$. 
Let $h\in D_{+}[a,b]$ or $h\in D_{-}[a,b]$ be increasing and strictly positive. Suppose that
$f^{(k)}(x)\geq (x-a)^{n-k-1}h(x)\frac{(n-1)!}{(n-k-1)!}$ for all $x\in [a,b]$. 
Then we have
\begin{equation*}
(b-a)f^{n+1}(a)g^n(a)h^{n-1}(a)
+\paren*{\int_a^bfgh}^n\leq \int_a^bf^{n+1}g^nh^{n-1}.
\end{equation*}
\end{proposition}

\begin{Proof}
We will only give a sketch of the proof. Assume that $h\in D_{+}[a,b]$ and define 
\begin{align*}
&F(x):=\paren*{\int_a^xfgh}^n-\int_a^xf^{n+1}g^nh^{n-1}, \\
&G(x):=n \paren*{\int_a^xfh}^{n-1}-f^n(x)h^{n-2}(x).
\end{align*}
Using the monotonicity of $g$, we see as before that $F^{\prime}\leq g^nfhG$.\par 
Next we define
\begin{equation*} 
H_s(x):=\paren*{\prod_{i=1}^s(n-i)}\paren*{\int_a^xfh}^{n-s-1}
-f^{n-s-1}(x)f^{(s)}(x)h^{n-s-2}(x)
\end{equation*}
for $s=1,\dots,k$.\par
Using similar arguments as before we obtain that $G_{+}^{\prime}\leq nfhH_1$
and $(H_{s})_{+}^{\prime}\leq fhH_{s+1}$ for all $s$.\par 
We further have, due to the monotonicity of $f$ and $h$,
\begin{equation*}
H_k(x)\leq f^{n-k-1}(x)h^{n-k-2}(x)\paren*{\frac{(n-1)!}{(n-k-1)!}(x-a)^{n-k-1}h(x)
-f^{(k)}(x)}.
\end{equation*}
Hence our assumption implies $H_k\leq 0$.\par  
It follows inductively that $H_s\leq 0$ for all $s$. Hence $G_{+}^{\prime}\leq 0$
and thus $G\leq G(a)=-f^n(a)h^{n-2}(a)$.\par 
This implies $F^{\prime}\leq -f^{n+1}(a)g^n(a)h^{n-2}(a)$ and the mean value theorem
gives us the desired conclusion.
\end{Proof}

\address
\email

\end{document}